\documentclass[12pt, leqno]{amsart}

\usepackage{amssymb,amsfonts,amscd,amsbsy}
\usepackage[mathscr]{eucal}

\theoremstyle{plain}
\newtheorem{thm}{Theorem}[section]
\newtheorem{lem}[thm]{Lemma}
\newtheorem{pro}[thm]{Proposition}

\theoremstyle{definition}
\newtheorem{Def}[thm]{Definition}
\newtheorem{rem}[thm]{Remark}

\begin{document}

\title[Densities of 4-ranks of $K_2(\mathcal{O})$]{Densities of 4-ranks of $K_2(\mathcal{O})$}

\author[Osburn]{Robert Osburn}

\address{Department of Mathematics, Louisiana State University,
Baton Rouge, LA 70803}

\email{osburn@math.lsu.edu}


\begin{abstract}
 In \cite{CH}, the authors established a method of
determining the structure of the 2-Sylow subgroup of the tame kernel $K_2(\mathcal{O})$ for certain
quadratic number fields. Specifically, the 4-rank for these fields was characterized in terms of
positive definite binary quadratic forms.  Numerical calculations led to questions concerning possible
density results of the 4-rank of tame kernels.  In this paper, we succeed in giving affirmative answers
to these questions.
\end{abstract}
\maketitle

\section{Introduction}
Since the $1960$'s, relationships between algebraic K-theory and number theory have been intensely
studied. For number fields F and their rings of integers $\mathcal{O}_{F}$, the K-groups
$K_{0}(\mathcal{O}_{F})$, $K_{1}(\mathcal{O}_{F})$, $K_{2}(\mathcal{O}_{F})$, $\dots$ were a main focus
of attention.  From \cite{Mil71} we have $$ K_{0}(\mathcal{O}_{F}) \cong \mathbb Z \times C(F) $$ where
$C(F)$ is the ideal class group of F, and $$ K_{1}(\mathcal{O}_{F}) \cong \mathcal{O}_{F}^{*}, $$ the
group of units of $\mathcal{O}_{F}$.

What can we say in general about $K_{2}(\mathcal{O}_{F})$?  Garland and Quillen in \cite{Gar71} and
\cite{Qu73} showed that $K_{2}(\mathcal{O}_{F})$ is finite.  A conjecture of Birch and Tate connects
the order of $K_{2}(\mathcal{O}_{F})$ and the value of the zeta function of F at $-1$ when F is a
totally real field.  For abelian number fields, this conjecture has been confirmed up to powers of 2
\cite{MW}. In \cite{tate} a 2-rank formula for $K_2(\mathcal{O}_F)$ was given by Tate.  Some results on
the 4-rank of $K_2(\mathcal{O}_F)$ were given in \cite{Qin1}, \cite{Qin2}, and \cite{Vaz}. To gain
further insight on the 4-rank of $K_2(\mathcal{O}_F)$, we consider the following specific families of
fields.

In this paper we deal with the 4-rank of the Milnor K-group $K_2(\mathcal{O})$ for the quadratic number
fields $\mathbb Q(\sqrt{pl})$, $\mathbb Q(\sqrt{2pl})$, $\mathbb Q(\sqrt{-pl})$, $\mathbb
Q(\sqrt{-2pl})$ for primes $p \equiv 7 \bmod 8$, $l \equiv 1 \bmod 8$ with $\Big(\frac{l}{p}\Big) =
1$.  In \cite{CH}, the authors show that for the fields $E = \mathbb Q(\sqrt{pl})$, $\mathbb
Q(\sqrt{2pl})$ and $F = \mathbb Q(\sqrt{-pl})$, $\mathbb Q(\sqrt{-2pl})$, $$ \mbox{4-rank}
\hspace{.05in} K_2(\mathcal{O}_E) = 1 \hspace{.05in} \mbox{or} \hspace{.05in} 2, $$ $$ \mbox{4-rank}
\hspace{.05in} K_2(\mathcal{O}_F) = 0 \hspace{.05in} \mbox{or} \hspace{.05in} 1. $$

Each of the possible values of 4-ranks is then characterized by checking which ones of the quadratic
forms $X^2 + 32Y^2$, $X^2 + 2pY^2$, $2X^2 + pY^2$ represent a certain power of $l$ over $\mathbb Z$.
This approach makes numerical computations accessible.  We should note that this approach involves
quadratic symbols and determining the matrix rank over $\mathbb F_2$ of $3 \times 3$ matrices with
Hilbert symbols as entries, see \cite{HK98}.  Fix a prime $p \equiv 7 \bmod 8$ and consider the set $$
\Omega = \{l \hspace{.05in}\mathrm{rational \hspace{.05in} prime}: l \equiv 1 \bmod 8 \hspace{.05in}
\mathrm{and} \hspace{.05in} \Big(\frac{l}{p}\Big) = \Big(\frac{p}{l}\Big) = 1 \}. $$

Let
\begin{center}

\item $\upsilon = $ $4$-rank $K_2(\mathcal{O}_{\mathbb Q(\sqrt{pl})})$

\item $\mu = $ $4$-rank $K_2(\mathcal{O}_{\mathbb Q(\sqrt{2pl})})$

\item $\sigma = $ $4$-rank $K_2(\mathcal{O}_{\mathbb Q(\sqrt{-pl})})$

\item $\tau = $ $4$-rank $K_2(\mathcal{O}_{\mathbb Q(\sqrt{-2pl})})$

\end{center}

and also consider the sets

\begin{center}

\item $\Omega_1 = \{l \in \Omega :  \upsilon = 1  \}$

\item $\Omega_2 = \{l \in \Omega :  \upsilon = 2 \}$

\item $\Omega_3 = \{l \in \Omega :  \mu = 1  \}$

\item $\Omega_4 = \{l \in \Omega :   \mu = 2 \}$

\item $\Lambda_1 = \{l \in \Omega :  \sigma = 0  \}$

\item $\Lambda_2 = \{l \in \Omega :  \sigma = 1 \}$

\item $\Lambda_3 = \{l \in \Omega : \tau  = 0  \}$

\item $\Lambda_4 = \{l \in \Omega :  \tau = 1 \}$.

\end{center}

We have computed the following (see Table 1 in Appendix): For $p=7$, there are 9730 primes $l$ in
$\Omega$ with $l \leq 10^6$. Among them, there are 4866 primes (50.01\%) in $\Omega_1$ and $\Omega_3$
and 4864 primes (49.99\%) in $\Omega_2$ and $\Omega_4$.  Also, there are 4878 primes (50.13\%) in
$\Lambda_1$ and $\Lambda_3$ and 4852 primes in $\Lambda_2$ and $\Lambda_4$. The goal of this paper is to prove the following theorem.

\begin{thm} \label{T:MQ0} For the fields $\mathbb Q(\sqrt{pl})$ and $\mathbb Q(\sqrt{2pl})$, 4-rank 1 and 2 each appear with natural density $\frac{1}{2}$ in
$\Omega$. For the fields $\mathbb Q(\sqrt{-pl})$ and $\mathbb Q(\sqrt{-2pl})$, 4-rank 0 and 1 each
appear with natural density $\frac{1}{2}$ in $\Omega$.
\end{thm}

Now consider the tuple of 4-ranks $(\upsilon$, $\mu$, $\sigma$, $\tau$).  By Corollary (5.6) in
\cite{CH}, there are eight possible tuples of 4-ranks.  For $p=7$, among the 9730 primes $l \in \Omega$
with $l \leq 10^6$, the eight possible tuples are realized by 1215, 1213, 1228, 1210, 1210, 1228, 1225,
1201 primes $l$ respectively (see Table 2 in Appendix). And, in fact:

\begin{thm} \label{T:MQ2} Each of the eight possible tuples of 4-ranks appear with natural density $\frac{1}{8}$ in $\Omega$.
\end{thm}

\section{Preliminaries}

Let $\mathcal{D}$ be a Galois extension of $\mathbb Q$, and $G = Gal(\mathcal{D}/ \mathbb Q)$.  Let
$Z(G)$ denote the center of G and $\mathcal{D}^{Z(G)}$ denote the fixed field of $Z(G)$.  Let $p$ be a
rational prime which is unramified in $\mathcal{D}$ and $\beta$ be a prime of $\mathcal{D}$ containing
$p$.  Let $\Big(\frac{\mathcal{D}/\mathbb Q}{p}\Big)$ denote the Artin symbol of $p$ and $\{g\}$ the
conjugacy class containing one element $g \in G$.

\begin{lem} \label{L:center} $\Big(\frac{\mathcal{D}/\mathbb Q}{p}\Big) = \{g\}$ for some $g \in Z(G)$ if and only if p splits completely
in $\mathcal{D}^{Z(G)}$.
\end{lem}

\begin{proof}
$\Big(\frac{\mathcal{D}/\mathbb Q}{p}\Big) = \{g\}$ for some $g \in Z(G)$ if and only if
$\Big(\frac{\mathcal{D}/\mathbb Q}{\beta}\Big) = g$ if and only if $\Big(\frac{\mathcal{D}^{Z(G)}/\mathbb Q}{\beta}\Big)$ = $\Big(\frac{\mathcal{D}/\mathbb Q}{\beta}\Big)\Big |_{\mathcal{D}^{Z(G)}}
= g |_{\mathcal{D}^{Z(G)}} = Id_{Gal(\mathcal{D}^{Z(G)}/ \mathbb Q)}$ if and only if p splits completely
in $\mathcal{D}^{Z(G)}$.
\end{proof}

Thus if we can show that rational primes split completely in the fixed field of the center of a certain
Galois group G, then we know the associated Artin symbol is a conjugacy class containing one element.
Hence we may identify the Artin symbol with this one element and consider the symbol to be an
automorphism which lies in Z(G).  Thus determining the order of Z(G) gives us the number of possible
choices for the Artin symbol.

Let $G_1$ and $G_2$ be finite groups and $A$ a finite abelian group. Suppose $r_1:G_1 \to A$ and
$r_2:G_2 \to A$ are two epimorphisms and $\mathcal{G} \subset G_1 \times G_2$ is the set $\{(g_1,
g_2)\in G_1 \times G_2 : r_1(g_1)=r_2(g_2) \}$.  Since A is abelian, there is an epimorphism $r:G_1
\times G_2 \to A$ given by $r(g_1, g_2)=r_1(g_1){r_2(g_2)}^{-1}$.  Thus $\mathcal{G} = ker(r) \subset
G_1 \times G_2$. One can check that $Z(\mathcal{G}) = \mathcal{G} \cap Z(G_1 \times G_2)$.

\begin{lem} \label{L:easy} (i) If $r_2\Big |_{Z(G_2)}$ is trivial, then $Z(\mathcal{G}) = ker(r_1\Big |_{Z(G_1)}) \times Z(G_2).$ \\
(ii) $Z(\mathcal{G}) = Z(G_1) \times Z(G_2) \iff r_1\Big |_{Z(G_1)}$ and $r_2\Big |_{Z(G_2)}$ are both trivial.
\end{lem}

\begin{proof}
(i) Suppose $(g_1, g_2) \in Z(\mathcal{G}) \subset ker(r)$ where $g_1 \in Z(G_1)$, $g_2 \in Z(G_2)$.
Thus $1 = r(g_1, g_2) = r_1(g_1){r_2(g_2)}^{-1}$ and so $r_1(g_1)=r_2(g_2)$. But $r_2(g_2) = 1$ which
yields $r_1(g_1) = 1$. Thus $g_1 \in ker(r_1\Big |_{Z(G_1)})$. The other inclusion is clear. \\
(ii) Take $(g_1, 1), (1, g_2) \in Z(G_1) \times Z(G_2) = Z(\mathcal{G}) \subset ker(r)$, respectively to
obtain that $r_1\Big |_{Z(G_1)}$ and $r_2\Big |_{Z(G_2)}$ are both trivial. The converse
follows from part (i).
\end{proof}

We will use the following definition throughout this paper.

\begin{Def} \label{D:sat} For primes $p \equiv 7\bmod 8$, $l \equiv 1\bmod 8$ with $\Big(\frac{l}{p}\Big) = \Big(\frac{p}{l}\Big)= 1$, $\mathcal{K} = \mathbb Q$$(\sqrt{-2p})$, and $h(\mathcal{K})$ the class number of $\mathcal{K}$, we say:

{\em l satisfies $\langle 1, 32\rangle$} if and only if $l = x^2 + 32y^2$ for some $x,y \in \mathbb Z$

{\em l satisfies $\langle 2, p\rangle$} if and only if $l^{\frac{h(\mathcal{K})}{4}} = 2n^{2} + pm^{2}$
for some $n,m \in \mathbb Z$ with $m \not \equiv 0\bmod l$

{\em l satisfies $\langle 1, 2p\rangle$} if and only if $l^{\frac{h(\mathcal{K})}{4}} = n^{2} + 2pm^{2}$ for some
$n,m \in \mathbb Z$ with $m \not \equiv 0\bmod l$.
\end{Def}

\section{Three Extensions}
In this section, we consider three degree eight field extensions of $\mathbb Q$.  The idea will be to
study composites of these fields and relate Artin symbols to 4-ranks.  Rational primes which split
completely in a degree 64 extension of $\mathbb Q$ will relate to Artin symbols and thus 4-ranks.
Therefore calculating the density of these primes will answer density questions involving 4-ranks.

\subsection{First Extension}

Consider $\mathbb Q(\sqrt{2})$ over $\mathbb Q$.  Let $\epsilon = 1 + \sqrt{2} \in (\mathbb
Z[\sqrt{2}])^{*}$. Then $\epsilon$ is a fundamental unit of $\mathbb Q(\sqrt{2})$ which has norm $-1$.
The degree 4 extension $\mathbb Q(\sqrt{2}, \sqrt{\epsilon}))$ over $\mathbb Q$ has normal closure
$\mathbb Q(\sqrt{2}, \sqrt{\epsilon}, \sqrt{-1})$.  Set $$ N_1 = \mathbb Q(\sqrt{2}, \sqrt{\epsilon},
\sqrt{-1}).$$ Note that $N_1$ is the splitting field of the polynomial $x^4-2x^2-1$ and so has degree 8
over $\mathbb Q$.  Therefore $Gal(N_1/\mathbb Q)$  is the dihedral group of order 8. Note that the automorphism induced by sending $\sqrt{\epsilon}$ to $-\sqrt{\epsilon}$ commutes with every element of
$Gal(N_1/\mathbb Q)$. Thus $Z(Gal(N_1/\mathbb Q)) = Gal(N_1/\mathbb Q(\sqrt{2}, \sqrt{-1}))$.

Observe that only the prime 2 ramifies in $\mathbb Q(\sqrt{2})$, $\mathbb Q(\sqrt{-1})$, $\mathbb
Q(\sqrt{\epsilon})$, and so only the prime 2 ramifies in the compositum $N_1$ over $\mathbb Q$. Now as
$l \in \Omega$ is unramified in $N_1$ over $\mathbb Q$, the Artin symbol $\Big(\frac{N_1/\mathbb Q}{\beta}\Big)$ is defined for primes $\beta$ of $\mathcal{O}_{N_1}$ containing $l$.  Let
$\Big(\frac{N_1/\mathbb Q}{l}\Big)$ denote the conjugacy class of $\Big(\frac{N_1/\mathbb Q}{\beta}\Big)$ in $Gal(N_1/\mathbb Q)$.  The primes ${l \in \Omega}$ split completely in $\mathbb Q(\sqrt{2},
\sqrt{-1})$ and ${N_1}^{Z(Gal(N_1/\mathbb Q))}  = \mathbb Q(\sqrt{2}, \sqrt{-1}).$ Thus by Lemma
\ref{L:center}, we have that $\Big(\frac{N_1/\mathbb Q}{l}\Big) = \{g\} \subset Z(Gal(N_1/\mathbb
Q))$. As $Z(Gal(N_1/\mathbb Q)$) has order 2, there are two possible choices for $\Big(\frac{N_1/\mathbb Q}{l}\Big)$.  Combining this statement with Addendum (3.4) from \cite{CH}, we have

\begin{rem} \label{R:N1}

\begin{eqnarray}
\Big(\frac{N_1/\mathbb Q}{l}\Big) = \{id\}
& \iff & l \hspace{.05in} \mbox{splits completely in} \hspace{.05in} N_1  \nonumber \\
& \iff & l \hspace{.05in} \mbox{satisfies} \hspace{.05in} \langle 1, 32\rangle. \nonumber
\end{eqnarray}

\end{rem}

\subsection{Second and Third Extension}
Consider the fixed prime $p \equiv 7\bmod 8$.  Note $p$ splits completely in $\mathcal{L} = \mathbb
Q(\sqrt{2})$ over $\mathbb Q$ and so $$ p\mathcal{O}_{\mathcal{L}} = \frak B \frak B^{'} $$ for some
primes $\frak B \not = \frak B^{'}$ in $\mathcal{L}$.  The field $\mathcal{L}$ has narrow class number
$h^{+}(\mathcal{L}) = 1$ as $h(\mathcal{L}) = 1$ and $N_{\mathcal{L}/\mathbb Q}(\epsilon) = -1$ where
$\epsilon = 1 + \sqrt{2}$ is a fundamental unit of $\mathbb Q(\sqrt{2})$, see \cite{Jan73}. From
\cite{CH},

\begin{lem} \label{L:pi} The prime $\frak B$ which occurs in the decomposition of $p\mathcal{O}_{\mathcal{L}}$ has a generator $\pi = a + b\sqrt{2} \in \mathcal{O}_{\mathcal{L}}$, unique up to a sign and to multiplication by the square of a unit in $\mathcal{O}_{\mathcal{L}}^{*}$ for which $N_{\mathcal{L}/\mathbb Q}(\pi) = a^2 - 2b^2 = -p$.
\end{lem}

The degree 4 extension $\mathbb Q(\sqrt{2}, \sqrt{\pi})$ over $\mathbb Q$ has normal closure

\noindent{$\mathbb Q(\sqrt{2}, \sqrt{\pi}, \sqrt{-p})$} as $N_{\mathcal{L}/\mathbb Q}(\pi) = -p$.  Set
$$ N_2 = \mathbb Q(\sqrt{2}, \sqrt{\pi}, \sqrt{-p}). $$ Then $N_2$ is Galois over $\mathbb Q$ and $[N_2
: \mathbb Q] = 8$.  Such an extension $N_2$ exists since the 2-Sylow subgroup of the ideal class group
of $\mathbb Q(\sqrt{-2p})$ is cyclic of order divisible by 4 \cite{CH3}.  Thus the Hilbert class field
of $\mathbb Q(\sqrt{-2p})$ contains a unique unramified cyclic degree 4 extension over $\mathbb
Q(\sqrt{-2p})$.  By Lemma 2.3 in \cite{CH}, $N_2$ is the unique unramified cyclic degree 4 extension
over $\mathbb Q(\sqrt{-2p})$.  Also compare \cite{JY87}. Similar to arguments in Section 2.1,
$Gal(N_2/\mathbb Q)$ is the dihedral group of order 8. Note that the automorphism induced by sending
$\sqrt{\pi}$ to $-\sqrt{\pi}$ commutes with every element of Gal$(N_2/\mathbb Q))$. Thus
$Z(Gal(N_2/\mathbb Q)) = Gal(N_2/\mathbb Q(\sqrt{2}, \sqrt{-p}))$.

\begin{pro} \label{P:ramN2} If $l \in \Omega$, then l is unramified in $N_2$ over $\mathbb Q$.
\end{pro}

\begin{proof} Since $p \equiv 7\bmod 8$, the discriminant of $\mathbb Q(\sqrt{-2p})$ is $-8p$.  For $l \in \Omega$, we have $\Big(\frac{-2p}{l}\Big) = 1$ and so $l$ is unramified in $\mathbb Q(\sqrt{-2p})$.  By Lemma 2.3 in \cite{CH}, we have $l$ is unramified in $N_2$ over $\mathbb Q$.
\end{proof}

As $l \in \Omega$ is unramified in $N_2$ over $\mathbb Q$, the Artin symbol $\Big(\frac{N_2/\mathbb Q}{\beta}\Big)$ is defined for primes $\beta$ of $\mathcal{O}_{N_2}$ containing $l$.  Let
$\Big(\frac{N_2/\mathbb Q}{l}\Big)$ denote the conjugacy class of $\Big(\frac{N_2/\mathbb Q}{\beta}\Big)$ in $Gal(N_2/\mathbb Q)$.  The primes $l \in \Omega$ split completely in $\mathbb Q(\sqrt{2},
\sqrt{-p})$ and ${N_2}^{Z(Gal(N_2/\mathbb Q))} = \mathbb Q(\sqrt{2}, \sqrt{-p}).$ By Lemma
\ref{L:center}, we have that $\Big(\frac{N_2/\mathbb Q}{l}\Big) = \{h\} \subset Z(Gal(N_2/\mathbb
Q))$ for some $h \in Z(Gal(N_2/\mathbb Q))$.  As Z(Gal$(N_2/\mathbb Q))$ has order 2, there are two
possible choices for $\Big(\frac{N_2/\mathbb Q}{l}\Big)$.  Combining this statement and Lemmas (3.3)
and (3.4) from \cite{CH}, we have

\begin{rem} \label{R:N2}

\begin{eqnarray}
\Big(\frac{N_2/\mathbb Q}{l}\Big) = \{id\}
& \iff & l \hspace{.05in} \mbox{splits completely in} \hspace{.05in} N_2 \nonumber \\
& \iff & l \hspace{.05in} \mbox{satisfies} \hspace{.05in} \langle 1, 2p\rangle. \nonumber
\end{eqnarray}

\begin{eqnarray}
\Big(\frac{N_2/\mathbb Q}{l}\Big) \not = \{id\}
& \iff & l \hspace{.05in} \mbox{does not split completely in} \hspace{.05in} N_2 \nonumber \\
& \iff & l \hspace{.05in} \mbox{satisfies} \hspace{.05in} \langle 2, p\rangle. \nonumber
\end{eqnarray}

\end{rem}

Finally, for $l \in \Omega$, $ l$ splits completely in $\mathbb Q(\zeta_{16})\iff l \equiv 1 \bmod 16$
This yields

\begin{rem} \label{R:cyc}

\begin{eqnarray}
\Big(\frac{\mathbb Q(\zeta_{16})/\mathbb Q}{l}\Big) = \{id\}
& \iff & l \hspace{.05in} \mbox{splits completely in} \hspace{.05in} \mathbb Q(\zeta_{16}) \nonumber \\
& \iff & l \equiv 1 \bmod 16. \nonumber
\end{eqnarray}

\end{rem}

\section{The Composite and Two Theorems}
In this section we consider the composite field $N_1N_2\mathbb Q(\zeta_{16})$. Set $$ L = N_1N_2\mathbb Q(\zeta_{16}). $$ Note that $[L:\mathbb Q] = 64$. As $N_1$, $N_2$, and $\mathbb
Q(\zeta_{16})$ are normal extensions of $\mathbb Q$, L is a normal extension of $\mathbb Q$.

For $l \in \Omega$, $l$ is unramified in $L$ as it is unramified in $N_1$, $N_2$, and $\mathbb
Q(\zeta_{16})$.  The Artin symbol $\Big(\frac{L/\mathbb Q}{\beta}\Big)$ is now defined for some
prime $\beta$ of $\mathcal{O}_L$ containing $l$.  Let $\Big(\frac{L/\mathbb Q}{l}\Big)$ denote the
conjugacy class of $\Big(\frac{L/\mathbb Q}{\beta}\Big)$ in $Gal(L/\mathbb Q)$.  Letting $M =
\mathbb Q(\sqrt{2}, \sqrt{-1}, \sqrt{-p}) \subset L$, we prove

\begin{lem} \label{L:Lcenter} $Z(Gal(L/\mathbb Q)) = Gal(L/M)$ is elementary abelian of order 8.
\end{lem}

\begin{proof} For $\sigma \in Gal(L/M)$, $\sigma$ can only change the sign of $\sqrt{\epsilon}$, $\sqrt{\pi}$, and $\sqrt{\zeta_{8}}$ as $\epsilon \in M$. Since $L=M(\sqrt{\epsilon}, \sqrt{\pi}, \sqrt{\zeta_{8}})$, $Gal(L/M)$ is elementary abelian of order
8. Now consider the restrictions $r_1:G_1 \to Gal(\mathbb Q(\sqrt{2})/\mathbb Q)$ and $r_2:G_2 \to
Gal(\mathbb Q(\sqrt{2})/\mathbb Q)$ where $G_1 = Gal(N_1/\mathbb Q)$ and $G_2 = Gal(N_2/\mathbb Q)$.
Clearly $r_1\Big |_{Z(G_1)}$ and

$r_1\Big |_{Z(G_2)}$ are both trivial.  Then by Lemma
\ref{L:easy} part (ii), $Z(\mathcal{G})$ is elementary abelian of order 4 where $\mathcal{G} =
Gal(N_1N_2/\mathbb Q)$. Now consider the restrictions $R_1:Gal(\mathbb Q(\zeta_{16})/\mathbb Q) \to
Gal(\mathbb Q(\zeta_{8})/\mathbb Q)$ and $R_2:\mathcal{G} \to Gal(\mathbb Q(\zeta_{8})/\mathbb Q)$.
Note that $ker(R_1)$ is cyclic of order 2 and $Z(\mathcal{G}) = Gal(M/\mathbb Q)$.  Thus $R_2\Big
|_{Z(\mathcal{G})}$ is trivial and so by the above and Lemma
\ref{L:easy} part (i), $Z(Gal(L/\mathbb Q)) \cong {\mathbb Z/2\mathbb Z} \times Z(\mathcal{G}) =
{\mathbb Z/2\mathbb Z} \times {\mathbb Z/2\mathbb Z} \times {\mathbb Z/2\mathbb Z}$. Thus
$Z(Gal(L/\mathbb Q)) = Gal(L/M)$.

\end{proof}

Now for $l \in \Omega$, $l$ splits completely in $\mathbb Q(\sqrt{-1})$ and $\mathbb Q(\sqrt{2},
\sqrt{-p})$ and so splits completely in the composite field $M = \mathbb Q(\sqrt{2}, \sqrt{-1},
\sqrt{-p})$.  From Lemma \ref{L:Lcenter}, $L^{Z(Gal(L/\mathbb Q))} = \mathbb Q(\sqrt{2}, \sqrt{-1},
\sqrt{-p}).$ So by Lemma \ref{L:center}, we have $\Big(\frac{L/\mathbb Q}{l}\Big) = \{k\} \subset
Z(Gal(L/\mathbb Q))$ for some $k \in Gal(L/\mathbb Q)$.  As $Z(Gal(L/\mathbb Q))$ has order 8, there
are eight possible choices for $\Big(\frac{L/\mathbb Q}{l}\Big)$.  Using Remarks \ref{R:N1},
\ref{R:N2}, and \ref{R:cyc}, we now make the following one to one correspondences.

\begin{rem} \label{R:eight}
(i) $\Big(\frac{L/\mathbb Q}{l}\Big) = \{id\} \iff$ $l$ splits completely in L  $\iff$ $\left \{
\begin{array}{l}
l \hspace{.05in} \mbox{splits completely in $N_1$,} \\
\mbox{$N_2$, and $\mathbb Q(\zeta_{16})$}
\end{array}
\right \}$ $\iff$ $\left \{ \begin{array}{l}
l \hspace{.05in} \mbox{satisfies} \hspace{.05in} \langle 1, 32\rangle \\
l \hspace{.05in} \mbox{satisfies} \hspace{.05in} \langle 1, 2p\rangle \\
l \equiv 1 \bmod 16
\end{array}
\right \} $.

(ii)  $\Big(\frac{L/\mathbb Q}{l}\Big) \not = \{id\} \iff$ $l$ does not split completely in L. Now
there are seven cases.

\begin{enumerate}
\item
$\left \{ \begin{array}{l}
l \hspace{.05in} \mbox{splits completely in $N_1$,} \\
\mbox{but does not in $N_2$ or $\mathbb Q(\zeta_{16})$}
\end{array}
\right \}$ $\iff$ $\left \{ \begin{array}{l}
l \hspace{.05in} \mbox{satisfies} \hspace{.05in} \langle 1, 32\rangle \\
l \hspace{.05in} \mbox{satisfies} \hspace{.05in} \langle 2, p\rangle \\
l \equiv 9 \bmod 16
\end{array}
\right \} $

\item
$\left \{ \begin{array}{l}
l \hspace{.05in} \mbox{splits completely in $N_1$} \\
\mbox{and $N_2$, but does not in $\mathbb Q(\zeta_{16})$}
\end{array}
\right \} $ $\iff$ $\left \{ \begin{array}{l}
l \hspace{.05in} \mbox{satisfies} \hspace{.05in} \langle 1, 32\rangle \\
l \hspace{.05in} \mbox{satisfies} \hspace{.05in} \langle 1, 2p\rangle \\
l \equiv 9 \bmod 16
\end{array}
\right \} $

\item
$\left \{ \begin{array}{l}
l \hspace{.05in} \mbox{splits completely in} \\
\mbox{$N_2$, but does not in $N_1$} \\
\mbox{or $\mathbb Q(\zeta_{16})$}
\end{array}
\right \} $ $\iff$ $\left \{ \begin{array}{l}
l \hspace{.05in} \mbox{does not satisfy} \hspace{.05in} \langle 1, 32\rangle \\
l \hspace{.05in} \mbox{satisfies} \hspace{.05in} \langle 1, 2p\rangle \\
l \equiv 9 \bmod 16
\end{array}
\right \}  $

\item
$\left \{ \begin{array}{l}
l \hspace{.05in} \mbox{splits completely in} \\
\mbox{$N_2$ and $\mathbb Q(\zeta_{16})$,} \\
\mbox{but does not in $N_1$}
\end{array}
\right \} $ $\iff$ $\left \{ \begin{array}{l}
l \hspace{.05in} \mbox{does not satisfy} \hspace{.05in} \langle 1, 32\rangle \\
l \hspace{.05in} \mbox{satisfies} \hspace{.05in} \langle 1, 2p\rangle \\
l \equiv 1 \bmod 16
\end{array}
\right \}  $

\item
$\left \{ \begin{array}{l}
l \hspace{.05in} \mbox{splits completely in $N_1$} \\
\mbox{and $\mathbb Q(\zeta_{16})$, but does not in $N_2$}
\end{array}
\right \} $ $\iff$ $\left \{ \begin{array}{l}
l \hspace{.05in} \mbox{satisfies} \hspace{.05in} \langle 1, 32\rangle \\
l \hspace{.05in} \mbox{satisfies} \hspace{.05in} \langle 2, p\rangle \\
l \equiv 1 \bmod 16
\end{array}
\right \}  $

\item
$\left \{ \begin{array}{l}
l \hspace{.05in} \mbox{splits completely in}  \\
\mbox{$\mathbb Q(\zeta_{16})$, but does not in $N_1$} \\
\mbox{or $N_2$}
\end{array}
\right \} $ $\iff$ $\left \{ \begin{array}{l}
l \hspace{.05in} \mbox{does not satisfy} \hspace{.05in} \langle 1, 32\rangle \\
l \hspace{.05in} \mbox{satisfies} \hspace{.05in} \langle 2, p\rangle \\
l \equiv 1 \bmod 16
\end{array}
\right \}  $

\item
$\left \{ \begin{array}{l}
l \hspace{.05in} \mbox{does not split completely} \\
\mbox{in $N_1$, $N_2$, or $\mathbb Q(\zeta_{16})$}
\end{array}
\right \} $ $\iff$ $\left \{ \begin{array}{l}
l \hspace{.05in} \mbox{does not satisfy} \hspace{.05in} \langle 1, 32\rangle \\
l \hspace{.05in} \mbox{satisfies} \hspace{.05in} \langle 2, p\rangle \\
l \equiv 9 \bmod 16
\end{array}
\right \}.  $

\end{enumerate}
\end{rem}

Now using Theorems (5.2), (5.3), (5.4), and (5.5) from \cite{CH}, we relate each Artin symbol
$\Big(\frac{L/\mathbb Q}{l}\Big)$ to each of the eight possible tuples of 4-ranks.

\begin{rem} \label{R:4ranks} {\em From Remark \ref{R:eight}, case (i) occurs if and only if we have $(2,2,1,1)$.  For case (ii), \\
(1) occurs if and only if we have $(1,2,0,1)$ \\
(2) occurs if and only if we have $(2,1,1,0)$ \\
(3) occurs if and only if we have $(2,1,0,1)$ \\
(4) occurs if and only if we have $(2,2,0,0)$ \\
(5) occurs if and only if we have $(1,1,0,0)$ \\
(6) occurs if and only if we have $(1,1,1,1)$ \\
(7) occurs if and only if we have $(1,2,1,0)$.

}
\end{rem}

We can now prove Theorem \ref{T:MQ2}.

\begin{proof}
Consider the set $ X = \{l \hspace{.05in} \mbox{prime}: l \hspace{.05in}\mbox{is unramified in $L$ and}
\Big(\frac{L/\mathbb Q}{l}\Big) = \{k\} \subset Z(Gal(L/\mathbb Q)) \} $ for some $k \in$
Gal$(L/\mathbb Q)$.  By the $\check C$ebotarev Density Theorem, the set $X$ has natural density
$\frac{1}{64}$ in the set of all primes. Recall $$ \Omega = \{l \hspace{.05in}\mathrm{rational
\hspace{.05in} prime}: l \equiv 1 \bmod 8 \hspace{.05in} \mathrm{and} \hspace{.05in} \Big(\frac{l}{p}\Big) = \Big(\frac{p}{l}\Big) = 1 \} $$ for some fixed prime $p \equiv 7 \bmod 8$.  By Dirichlet's Theorem on primes in arithmetic progressions, $\Omega$ has natural density $\frac{1}{8}$ in the set of
all primes.  Thus $X$ has natural density $\frac{1}{8}$ in $\Omega$.  By Remark \ref{R:eight} and
\ref{R:4ranks}, each of the eight choices for $\Big(\frac{L/\mathbb Q}{l}\Big)$ is in one to one
correspondence with each of the possible tuples of 4-ranks.  Thus each of the eight possible tuples of
4-ranks appear with natural density $\frac{1}{8}$ in $\Omega$.

\end{proof}

Now we can prove Theorem \ref{T:MQ0}

\begin{proof}
We see from Remark \ref{R:4ranks}, $4$-rank $K_2(\mathcal{O}_{\mathbb Q(\sqrt{pl})}) = 1$ in cases (ii),
parts (1), (5), (6), and (7), $4$-rank $K_2(\mathcal{O}_{\mathbb Q(\sqrt{2pl})}) = 2$ in case (i) and
case (ii) parts (1), (4), and (7), $4$-rank $K_2(\mathcal{O}_{\mathbb Q(\sqrt{-pl})}) = 0$ in case (ii)
parts (1), (3), (4), and (5), $4$-rank $K_2(\mathcal{O}_{\mathbb Q(\sqrt{-2pl})}) = 1$ in case (i) and
case (ii) parts (1), (3), and (6).  As each of the 4-rank tuples occur with natural density
$\frac{1}{8}$, we have for the fields $\mathbb Q(\sqrt{pl})$ and $\mathbb Q(\sqrt{2pl})$, 4-rank 1 and
2 each appear with natural density $4\cdot\frac{1}{8} = \frac{1}{2}$ in $\Omega$. For the fields
$\mathbb Q(\sqrt{-pl})$ and $\mathbb Q(\sqrt{-2pl})$, 4-rank 0 and 1 each appear with natural density
$4\cdot\frac{1}{8} = \frac{1}{2}$ in $\Omega$.

\end{proof}

\section*{Appendix} The following tables motivated possible
density results of 4-ranks of tame kernels.  We consider primes $l \in \Omega$ with $l \leq N$ for a
fixed prime $p \equiv 7 \bmod 8$ and positive integer N.  For Table 1, we consider the sets $\Omega_1,
\dots ,\Omega_4$ and $\Lambda_1, \dots ,\Lambda_4$ as in the Introduction.  For Table 2, we consider the
sets

\begin{center}

\item $I_1 = \{l \in \Omega :$ 4-rank tuple is (1,1,0,0)\}

\item $I_2 = \{l \in \Omega :$ 4-rank tuple is (1,1,1,1)\}

\item $I_3 = \{l \in \Omega :$ 4-rank tuple is (2,1,1,0)\}

\item $I_4 = \{l \in \Omega :$ 4-rank tuple is (2,1,0,1)\}

\item $I_5 = \{l \in \Omega :$ 4-rank tuple is (1,2,1,0)\}

\item $I_6 = \{l \in \Omega :$ 4-rank tuple is (1,2,0,1)\}

\item $I_7 = \{l \in \Omega :$ 4-rank tuple is (2,2,0,0)\}

\item $I_8 = \{l \in \Omega :$ 4-rank tuple is (2,2,1,1)\}.

\end{center}

\section*{Table 1}
\begin{center}
\begin{tabular}{|c||c|c|c|c|c|c|}
\hline Primes & $p=7$ & & $p=23$ & & $p=31$ &  \cr \hline \hline Cardinality & $N = 1000000$ & $\%$ &
$N = 1000000$ &$\%$ & $N = 1000000$ & $\%$  \cr \hline \hline $\mid \Omega \mid$   & 9730 & &  9742 & &
9754 &   \cr \hline $\mid \Omega_1 \mid$ & 4866 & 50.01 & 4905 & 50.35 & 4916 & 50.40  \cr
   \hline
$\mid \Omega_2 \mid$ & 4864 & 49.99 & 4837 & 49.65 & 4838 & 49.60 \cr
   \hline
$\mid \Omega_3 \mid$ & 4866 & 50.01 & 4911 & 50.41 & 4851 & 49.73 \cr
   \hline
$\mid \Omega_4 \mid$ & 4864 & 49.99 & 4831 & 49.59 & 4903 & 50.27 \cr \hline $\mid \Lambda_1 \mid$ &
4878 & 50.13 & 4912 & 50.42 & 4930 & 50.54  \cr \hline $\mid \Lambda_2 \mid$ & 4852 & 49.87 & 4830 &
49.58 & 4824 & 49.46  \cr \hline $\mid \Lambda_3 \mid$ & 4878 & 50.13 & 4876 & 50.05 & 4943 & 50.68 \cr
\hline $\mid \Lambda_4 \mid$ & 4852 & 49.87 & 4866 & 49.95 & 4811 & 49.32 \cr \hline

\end{tabular}
\end{center}

\section*{Table 2}
\begin{center}
\begin{tabular}{|c||c|c|c|c|c|c|}
\hline Primes & $p=7$ & & $p=23$ & & $p=31$ &  \cr \hline \hline Cardinality & $N = 1000000$ & $\%$ &
$N = 1000000$ &$\%$ & $N = 1000000$ & $\%$  \cr \hline \hline
   $\mid \Omega \mid$ & 9730 & & 9742 & & 9754 &  \cr \hline
   $\mid I_1 \mid$ & 1215 & 12.49 & 1246 & 12.79 & 1246 & 12.77 \cr
   \hline
   $\mid I_2 \mid$ & 1213 & 12.46 & 1229 & 12.62 & 1203 & 12.33 \cr
   \hline
   $\mid I_3 \mid$ & 1228 & 12.62 & 1211 & 12.43 & 1214 & 12.45 \cr
   \hline
   $\mid I_4 \mid$ & 1210 & 12.44 & 1225 & 12.57 & 1188 & 12.18 \cr
   \hline
   $\mid I_5 \mid$ & 1210 & 12.44 & 1204 & 12.36 & 1227 & 12.58 \cr
   \hline
   $\mid I_6 \mid$ & 1228 & 12.62 & 1226 & 12.58 & 1240 & 12.71 \cr
   \hline
   $\mid I_7 \mid$ & 1225 & 12.59 & 1215 & 12.47 & 1256 & 12.88 \cr
   \hline
   $\mid I_8 \mid$ & 1201 & 12.34 & 1186 & 12.17 & 1180 & 12.10 \cr \hline

\end{tabular}
\end{center}

\section*{acknowledgments} I would like to thank the referee for the valuable comments and suggestions.  I also deeply thank P.E. Conner and my advisor, J. Hurrelbrink.


\begin{thebibliography}{10}


\bibitem{CH}
P. E. Conner and J. Hurrelbrink, \emph{On the 4-rank of the tame kernel $K_2(\mathcal{O})$ in positive
definite terms}, J. Number Th. \textbf{88} (2001), 263--282.

\bibitem{CH3}
P. E. Conner and J. Hurrelbrink, {\em Class Number Parity}, Ser. Pure Math. \textbf{8}, World Sci.
Publ. Singapore, 1988.

\bibitem{Gar71}
H. Garland, \emph{A finiteness theorem for $K_2$ of a number field}, Ann. of Math. \textbf{94} (1971),
534--548.

\bibitem{HK98}
J. Hurrelbrink and M. Kolster, \emph{Tame kernels under relative quadratic extensions and Hilbert
symbols}, J. reine angew. Math. \textbf{499} (1998), 145--188.

\bibitem{Jan73}
G. Janusz, {\em Algebraic Number Fields}, Academic Press, New York, 1973.

\bibitem{JY87}
C. U. Jensen and N. Yui, \emph{Quaternion Extensions}, Algebraic Geometry and Commutative Algebra in
Honor of Masayoshi Nagata, 1987, 155--182.

\bibitem{MW}
B. Mazur and A. Wiles, \emph{Class fields of abelian extensions of $\mathbb Q$}, Invent. Math.
\textbf{76} (1984), 179--330.

\bibitem{Mil71}
J. Milnor, {\em An Introduction to Algebraic K-Theory}, Ann. Math. Studies. Vol. \textbf{72}, Princeton
Univ. Press, Princeton, 1971.

\bibitem{Qin1}
H. Qin, \emph{The 2-Sylow subgroups of the tame kernel of imaginary quadratic fields}, Acta Arith.,
\textbf{69} (1995), 153--169.

\bibitem{Qin2}
H. Qin, \emph{The 4-rank of $K_2(\mathcal{O}_F)$ for real quadratic fields F}, Acta Arith., \textbf{72}
(1995), 323--333.

\bibitem{Qu73}
D. Quillen, \emph{Higher algebraic K-theory}, Algebraic K-Theory I, Lecture Notes in Mathematics, Vol.
\textbf{341}, Springer-Verlag, New York, 1973, 85--147.

\bibitem{tate}
J. Tate, \emph{Relations between $K_2$ and Galois cohomology}, Inventiones Math., \textbf{36}, (1976),
257--274.

\bibitem{Vaz}
A. Vazzana, \emph{4-ranks of $K_2$ of rings of integers in quadratic number fields}, Ph.D. Thesis,
University of Michigan, 1998.

\end{thebibliography}
\end{document}